\documentclass[a4paper, 11pt, reqno]{amsart}

\usepackage{amsthm,amssymb,amstext,amscd,amsfonts,amsbsy,amsrefs,amsxtra,latexsym,amsmath,xcolor,mathrsfs,fancybox,upgreek,soul,url}
\usepackage[english]{babel}
\usepackage[all,cmtip]{xy}
\usepackage[latin1]{inputenc}
\usepackage{hyperref}
\usepackage{comment}
\usepackage{mdframed}
\allowdisplaybreaks
\usepackage{mathtools}
\usepackage{enumerate}
\usepackage{longtable}
\usepackage{thmtools}
\usepackage{thm-restate}
\usepackage{chngcntr}
\usepackage{enumitem}
\usepackage{etoolbox}
\usepackage{todonotes}
\usepackage{nicefrac}

\usepackage{graphicx}

\usepackage{tikz}
\usepackage{verbatim}
\usetikzlibrary{quotes,angles}

\DeclarePairedDelimiter\abs{\lvert}{\rvert}
\DeclarePairedDelimiter\norm{\lVert}{\rVert}

\makeatletter
\let\oldabs\abs
\def\abs{\@ifstar{\oldabs}{\oldabs*}}
\let\oldnorm\norm
\def\norm{\@ifstar{\oldnorm}{\oldnorm*}}
\makeatother

\makeatletter
\glb@settings 
\makeatother

\newtheorem{theorem}{Theorem}

\newtheorem{corollary}[theorem]{Corollary}

\newtheorem{conjecture}[theorem]{Conjecture}

\theoremstyle{definition}

\theoremstyle{remark}

\newtheorem*{remark}{Remark}

\numberwithin{theorem}{section}
\numberwithin{proposition}{section}
\numberwithin{lemma}{section}
\numberwithin{corollary}{section}
\numberwithin{equation}{section}
\numberwithin{conjecture}{section}
\numberwithin{definition}{section}

\setlist[enumerate,1]{before=}
\AfterEndEnvironment{enumerate}{}

\newcommand{\N}{\mathbb{N}}
\newcommand{\Z}{\mathbb{Z}}

\usepackage{bm}

\author{Johann Stumpenhusen}
\title{On the solutions of $\varphi(dn) = \varphi(d(n + h))$}
\address{Department of Mathematics and Computer Science, Division of Mathematics, University of Cologne, Weyertal 86-90, 50931 Cologne, Germany}
\email{jstumpen@math.uni-koeln.de}
\date{\today}

\begin{document}

\maketitle

\begin{abstract}
    This is a short note about a chapter in the author's bachelor thesis regarding a paper by Ford \cite{Ford} concerning a conjecture by Erd\H{o}s.
\end{abstract}

\section{Introduction and Statement of Results}

Let $\varphi$ be Euler's totient function and $\sigma$ the sum-of-divisors-function. Ford \cite{Ford} picked up an unsolved problem which is based on a conjecture by Erd\H{o}s who asked whether there are infinitely many $n \in \N$ satisfying the equations
\[\varphi(n) = \varphi(n + 1)\]
or
\[\sigma(n) = \sigma(n + 1).\]
These problems have been treated in a variety of preceding papers, often adding an arbitrary but fixed $k \in \N$ instead of 1. For this, see the references in Ford's article. The main idea presented in \cite{Ford} is to use a celebrated theorem by Maynard, approaching Schinzel's hypothesis $H$.

\begin{conjecture}[Schinzel, Sierpi\'nski \cite{Schinzel}] Let $m \in \N$ and $f_1, \ldots, f_m \in \Z[x]$ be a family of polynomials. If there is no prime number $p \in \mathbb{P}$ such that $p \mid \prod_{i = 1}^m f_i(n)$ for all $n \in \N$, then there are infinitely many $n \in \N$ such that the numbers $f_1(n), \ldots, f_m(n)$ all are simultaneously prime.
\end{conjecture}

If a set of polynomials satisfies the condition of the conjecture, it is said to be \textit{admissible}. The special case of $m = 2, f_1(x) = x$ and $f_2(x) = x + 2$ is also known as the Twin Prime Conjecture. Maynard proved the following weakened version of this conjecture.

\begin{theorem}[Maynard {\cite[Theorem 1.2]{Maynard15}}] \label{thm:Maynard15} Let $m \in \N$. Then there exists $k_m \in \N$ such that for every admissible set $f_1, \ldots, f_{k_m} \in \Z[x]$ of $k_m$ linear polynomials, there are infinitely many $n \in \N$ such that at least $m$ of the numbers $f_1(n), \ldots, f_{k_m}(n)$ are simultaneously prime.
\end{theorem}

Using a result by the PolyMath8b project \cite{PolyMath8b}, Ford used the fact that $k_2 \leq 50$ to prove the following theorem.

\begin{theorem}[Ford {\cite[Theorem 1]{Ford}}]\label{thm:Ford}  Let $l \in \N$ be arbitrary.
    \begin{itemize}
        \item[(a)] Let $\Lambda := 2^53^35^2\prod_{7 \leq p \leq 47}p$. Then there exist infinitely many $n \in \N$ satisfying
        \[\sigma(n) = \sigma(n + l\Lambda).\]
        \item[(b)] There exists an even $h \leq 3570$ such that the equation
        \[\sigma(n) = \sigma(n + lh)\]
        is satisfied for infinitely many $n \in \N$.
    \end{itemize}
\end{theorem}

We will tweak the set used in Ford's proof and adjust the rest of it a bit as well to derive the following result.

\begin{theorem}\label{thm:phi(d)phi(n)}
    Let $d \in \N$ and $\kappa_d$ be given by
    \begin{align*}
    \kappa_d &:= \max\left\{d\frac{k_1 - k_2}{(k_1d + 1,k_2d + 1)}\prod_{p \mid \frac{(k_1d + 1)(k_2d + 1)}{(k_1d + 1, k_2d + 1)}}p:0 \leq k_1,k_2 \leq 49\right\}\\
    &< 49d(48d+1)(49d+1) \ll d^3.
    \end{align*}
    There exists an $h \leq \kappa_d$ with $d \mid h$ such that for every $l \in \N$ with $(d,l) = 1$ there exist infinitely many $n \in \N$ satisfying
    \[\varphi(d)\varphi(n) = \varphi(dn) = \varphi(d(n + lh)).\]
    
\end{theorem}

\begin{corollary}\label{cor:phi(n)=phi(2n+2lh)}
$\kappa_d$ takes the following first few values:
\begin{center}
\begin{tabular}{c|c}
    $d$ & $\kappa_d$ \\
    \hline
    2 & 227'950\\
    3 & 762'120\\
    4 & 1'910'900\\
    5 & 3'705'990\\
    6 & 6'414'480\\
    7 & 9'506'770\\
    8 & 15'169'800\\
    9 & 21'580'650\\
    10 & 29'582'750\\
    11 & 34'729'310\\
    12 & 51'059'232
\end{tabular}
\begin{tabular}{c|c}
    13 & 64'900'550\\
    14 & 81'031'650\\
    15 & 93'549'750\\
    16 & 12'0877'440\\
    17 & 144'970'050\\
    18 & 172'052'550\\
    19 & 190'261'896\\
    20 & 208'224'480\\
    21 & 272'127'030\\
    22 & 313'911'312\\
    23 & 337'307'880\\
    24 & 407'477'400\\
    25 & 460'516'250
\end{tabular}
\begin{tabular}{c|c}
    26 & 486'298'150\\
    27 & 545'455'944\\
    28 & 646'820'300\\
    29 & 675'761'016\\
    30 & 795'443'250\\
    31 & 825'322'920\\
    32 & 965'248'800\\
    33 & 1'058'536'050\\
    34 & 1'157'648'150\\
    35 & 1'110'432'750\\
    36 & 1'373'988'960\\
    37 & 1'491'697'550\\
    38 & 1'425'958'094
\end{tabular}
\begin{tabular}{c|c}
    39 & 1'642'680'936\\
    40 & 1'884'442'560\\
    41 & 2'022'423'810\\
    42 & 2'181'407'550\\
    43 & 2'201'405'640\\
    44 & 2'507'943'900\\
    45 & 2'682'771'750\\
    46 & 2'865'437'520\\
    47 & 2'874'316'584\\
    48 & 3'255'610'800\\
    49 & 3'463'263'650\\
    50 & 3'679'563'750\\
    51 & 3'665'785'650
    \end{tabular}
\end{center}
\end{corollary}

\section{Proofs}

\begin{proof}[Proof of Theorem \ref{thm:phi(d)phi(n)}] Let
\[A(d) := \{kd + 1: k \in \N, 0 \leq k \leq 49\}.\]
The set $\mathcal{L} = \{ax + 1: a \in A(d)\} \subset \Z[x]$ is admissible. To see this, set
\[P(x) := \prod_{k = 0}^{49} \left((kd + 1) \cdot x + 1\right)\]
Thus $\left(P(1),P(P(1))\right) = 1$ by construction and hence $\mathcal{L}$ is admissible.

Similarly to the proof of Theorem \ref{thm:Ford} in \cite{Ford}, we deduce from Theorem \ref{thm:Maynard15} that there are at least two elements $a_1 > a_2 \in A$ such that $a_1n + 1$ and $a_2n + 1$ are simultaneously prime infinitely often, directly implying that the same holds for $a_1'n + 1$ and $a_2'n + 1$ where $a_i' := \frac{a_i}{(a_1,a_2)}$ for $i \in \{1,2\}$ via multiplying $n$ by $(a_1,a_2)$. We then set
\[m_1 := a_2'ls(a_1,a_2)(a_1'r + 1), \quad m_2 := a_1'ls(a_1,a_2)(a_2'r + 1)\]
 where $s(a_1,a_2) := \prod_{p \mid a_1'a_2'}p$ and $r$ is chosen in a way such that $a_1'r + 1$ and $a_2'r + 1$ are prime and $r > \max\{a_1,l\}$. Note that all these factors are coprime to $d$ by construction.  Then we have
\[m_1 - m_2 = ls(a_1,a_2)(a_1' - a_2').\]
In particular, $d \mid a_1' - a_2' = \frac{k_1d + 1 - k_2d - 1}{(k_1d + 1, k_2d + 1)} = \frac{k_1 - k_2}{(k_1d + 1, k_2d + 1)}d$ by construction\footnote{Note that $d$ is coprime to the denominator.} for some $0 \leq k_1, k_2 \leq 49$, so $d \mid m_1 - m_2$. Furthermore, we have
\begin{align*}
    \varphi(dm_1) &= \varphi(d)\varphi\left(a_2'ls(a_1,a_2)(a_1'r + 1)\right)\\
    &= \varphi(d)a_2'\varphi(ls(a_1,a_2))\varphi(a_1'r + 1)\\
    &= \varphi(d)a_1'a_2'r\varphi(ls(a_1,a_2))\\
    &= \varphi(dm_2),
\end{align*}
yielding the claim. For the bound on $\kappa_d$, we see that we just need to calculate the greatest possible value of
\[\frac{m_1 - m_2}{l}\]
while varying $a_1, a_2 \in A(d)$. Using the most trivial bounds on each factor of the product in the set of Theorem \ref{thm:phi(d)phi(n)}, we retrieve the presented upper bound.\footnote{A java code may be shared on request.}
\end{proof}

\begin{proof}[Proof of Corollary \ref{cor:phi(n)=phi(2n+2lh)}]
    We simply evaluate the maximum of the set given in Theorem \ref{thm:phi(d)phi(n)}.
\end{proof}

\begin{remark}
    After an explosive growth for small $d$, it appears that our bound for $\kappa_d$ increases steadily slower as $d \to \infty$. More interestingly, our bound is not monotone.
\end{remark}

\end{document}